\documentclass{article}

\usepackage{amsfonts}

\usepackage{amsmath}
\usepackage{amsthm}
\usepackage{amssymb}

\providecommand\mathbb[1]{\mathsf{##1}}
\providecommand\mathfrak[1]{\mathcal{##1}}

\newtheorem{thm}{Theorem}[section]
\newtheorem{lem}[thm]{Lemma} 
\newtheorem{prop}[thm]{Proposition} 
\newtheorem{cor}[thm]{Corollary} 

\newcommand{\BZ}{{\mathbf{Z}}}
\newcommand{\BQ}{{\mathbf{Q}}}

\newcommand{\abs}[1]{\left| #1 \right|}

\newcommand{\OK}{{\mathcal K}}

\newcommand{\OO}{{\mathcal O}}
\newcommand{\ga}{{\mathfrak a}}

\newcommand{\gm}{{\mathfrak m}}
\newcommand{\gp}{{\mathfrak p}}

\newcommand{\floor}[1]{\left\lfloor#1\right\rfloor}

\begin{document}

\title{Multiplicative structures of values of the sum-of-divisors function
}

\author{Tomohiro Yamada\\
Department of Mathematics,\\
University of Kyoto, Kyoto, Japan}

\maketitle

\begin{abstract}
We study equations of the form $\sigma(p^{q-1})=Az$, where $p$ is a prime, $q$ is a fixed odd prime, $A$ is a fixed integer and $z$ is an integer composed of primes in a fixed finite set.
We shall improve upper bounds for the size and the number of solutions of such equations.
\end{abstract}

\section{Introduction}\label{intro}
We denote by $\sigma(N)$ the sum of divisors of $N$ a positive integer.
Then $\sigma(N)$ is multiplicative and $\sigma(p^{a-1})=(p^a-1)/(p-1)$
for any positive prime $p$ and positive integer $a$.  Thus we are led to
study numbers of the form $(p^a-1)/(p-1)$.

In this paper, we first study the equation
\begin{equation}\label{eq11}
\frac{x^q-1}{x-1}=Am_1^{e_1}m_2^{e_2}\cdots m_s^{e_s},
\end{equation}
where $x, m_1, \cdots, m_s$ are positive rational (not necessarily prime) integers, $q$
is a positive prime.  We note that if $x$ is prime, then the left of (\ref{eq11}) is equal to
$\sigma(x^q)$.  A considerable result in this direction is \cite[Theorem 5]{BHM},
which states that if (\ref{eq11}) with $A=s=1$ and $e_1$ prime holds, then $e_1\leq 9000q^2\log^4 q$.

We use a similar argument to \cite{BHM} to obtain our main theorem,
which improves the upper bound in \cite{BG2}.  Before stating this result, 
we introduce some notations.  For real $x$, we denote by $\langle x\rangle$
the quantity $\max\{x, 2\}$.  Moreover, let $c_2(m)=1500\cdot 38^{m+1}(m+1)^{3m+9}$.

\begin{thm}\label{thm1}
Let $A$ be a positive integer, $q$ be a positive prime, and
$m_1<m_2<\cdots<m_s$ be positive integers which are composed by primes congruent to $1$ mod $q$.
Denote by $h$ the class number of the quadratic field $\BQ(\sqrt{(-1)^{(q-1)/2}q})$.
Moreover, in the case $q\equiv 1\pmod{4}$, denote by $R$ the regulator of this field.
If a positive integer $x$ satisfies the equation (\ref{eq11}), then $e_i<U$, where
we denote by $U$ the constant
\begin{equation}\label{ineq11a}
\max\{\frac{qh^2\log(4m_s^{1/2})}{\log m_1}, \frac{q\log 2q^{1/2}}{\log m_1}, C_0\langle\log A^\prime+2R\rangle(\log\frac{2s(s+2)C_0}{h})\}
\end{equation}
in the case $q\equiv 1\pmod {4}$, and
\begin{equation}\label{ineq11b}
\max\{\frac{q\log 2q^{1/2}}{\log m_1}, C_1\langle\log A^\prime\rangle(\log\frac{(s+1)C_1}{h})\}+h-1
\end{equation}
in the case $q\equiv 3\pmod {4}$, where
\begin{equation}
C_0=2^4qc_2(s+2)\langle R\rangle(1+\frac{2R}{\log m_1})\prod_{i=2}^{s}(\log m_i+2R),
\end{equation}
\begin{equation}
C_1=\frac{2^6qc_2(s+1)}{\log 7}\prod_{i=2}^{s}(\log m_i),
\end{equation}
and $A^\prime=A m_1^{v_1}m_2^{v_2}\cdots m_s^{v_s}$ for some $v_1, v_2, \cdots, v_s\leq h-1$.
\end{thm}

We derive this theorem from the following theorem concerning values of binary quadratic forms.

\begin{thm}\label{thm2}
Let $A$ be a positive integer, $D$ be an integer with $D\equiv 1 \pmod {4}$, and
$m_1<m_2< \cdots<m_s$ be positive integers which are composed by primes congruent to $1$ mod $D$.
Let $c_0$ and $c_1$ be positive real numbers with $c_1<1/2$.
Denote by $h$ the class number of the quadratic field $\BQ((-1)^{(D-1)/2}\sqrt{D})$.
Moreover, in the case $q\equiv 1\pmod{4}$, denote by $R$ the regulator of this field.
If $X, Y$ are two integers satisfying

\begin{equation}\label{eq12}
X^2-DY^2=4Am_1^{e_1}m_2^{e_2}\cdots m_s^{e_s},
\end{equation}
\begin{equation}\label{ineq12}
\abs{Y}<c_0\abs{Am_1^{e_1}m_2^{e_2}\cdots m_s^{e_s}}^{1/2-c_1},
\end{equation}
and any prime ideal dividing $(\frac{X+Y\sqrt{D}}{2}, \frac{X-Y\sqrt{D}}{2})$
divides also $\sqrt{D}$, then

\begin{equation}\label{ineq13a}
e_i<\max\{\frac{h^2\log(2c_0m_s^{1/2})}{c_1\log m_1}, \frac{\log c_0\abs{D}^{1/2}}{c_1\log m_1}, C_2\langle\log A^\prime+2R\rangle(\log\frac{2s(s+2)C_2}{h})\}
\end{equation}
in the case $D>0$, and
\begin{equation}\label{ineq13b}
e_i<\max\{\frac{\log c_0\abs{D}^{1/2}}{c_1\log m_1}, C_3\langle\log A^\prime\rangle(\log\frac{(s+1)C_3}{h})\}+h-1
\end{equation}
in the case $D<0$, where
\begin{equation}
C_2=\frac{2^4 c_2(s+2)}{c_1}\langle R\rangle(1+\frac{2R}{\log m_1})\prod_{i=2}^{s}(\log m_i+2R),
\end{equation}
\begin{equation}
C_3=\frac{2^6 c_2(s+1)}{c_1\log 7}\prod_{i=2}^{s}(\log m_i),
\end{equation}
and $A^\prime=A m_1^{v_1}m_2^{v_2}\cdots m_s^{v_s}$ for some $v_1, v_2, \cdots, v_s\leq h-1$.
\end{thm}

We next consider the largest prime divisor of the left of (\ref{eq11}).  Denote by $P[n]$
the largest prime divisor of an integer $n$.  Kotov\cite{Kot} shows that $P[(x^q-1)/(x-1)]>c\log \log x$ for some effectively
computable constant $c>0$ depending only on $q$.  An explicit estimate can be found
in \cite[Theorem 3]{Gyo}.  See also a general result of \cite{BG2}.

Using Theorem \ref{thm1}, we can improve known results on the largest prime divisor of
the left of (\ref{eq11}).

\begin{thm}\label{thm3}
Let $x$ be a positive integers and $q$ be a positive prime.  Then, for any real $\epsilon>0$,
there exists an effectively computable constant $x_0$ depending only on $q$ and $\epsilon$
such that
\begin{equation}\label{eq13}
P[\frac{x^q-1}{x-1}]>(\frac{q-1}{6}-\epsilon)\log\log x
\end{equation}
for any integer $x>x_0$.
\end{thm}

Finally, we consider the number of solutions of the equation
\begin{equation}\label{eq14}
\frac{p^q-1}{p-1}=m_1^{e_1}m_2^{e_2}\cdots m_s^{e_s},
\end{equation}
where $p, q, m_1, \cdots, m_s$ are positive rational integers with $p$ and $q$ prime.
We note that the left of this equation is equal to $\sigma(p^{q-1})$.

\begin{thm}\label{thm4}
Set $c_7=2^{12}\times 38^2\times 1500^2$.  If $q>\frac{16}{9}es^4$ is prime,
then the equation (\ref{eq14}) has at most
\begin{equation}\label{eq14b}
s(\frac{\log c_7+19s\log(s+2)+3\sum_{i=2}^{s}\log\log m_i}{\log q}+7)
\end{equation}
solutions in integers $p, e_1, e_2, \cdots, e_s$ with $p$ prime.
\end{thm}

\section{lower bounds for linear forms in logarithms}
Our main tool is lower bounds for linear forms in logarithms of some special form.
We denote by $H(\alpha)$ the absolute height of $\alpha$(See Section 5 in \cite{BG1}) and set $h(\alpha)=\log H(\alpha)$.
Let $\OK$ be an algebraic field of degree $k$ over $\BQ$ and $\alpha_1, \cdots, \alpha_m(m\ge 2)$ be nonzero algebraic numbers in $\OK$.
Let $h_1, \cdots, h_m$ be real numbers such that
\begin{equation}\label{ineq21}
h_i\ge\max\{h(\alpha_i), \abs{\frac{\log\alpha_i}{3.3k}}, \frac{1}{k}\} \mbox{ for } i=1, \cdots, m,
\end{equation}
where log denotes the principal value of the logarithm.  Let $b_1, \cdots, b_{m-1}$ be
rational integers and put $B=\max\{\abs{b_1}, \cdots, \abs{b_{m-1}}, 3\}.$  Set

\begin{equation*}
\Lambda=\alpha_1^{b_1}\cdots\alpha_{m-1}^{b_{m-1}}\alpha_m-1.
\end{equation*}
The following estimate is due to \cite{BG1} and \cite{Wal}.
\begin{prop}\label{prop1}
If $\lambda\ne 0$,
\begin{equation}\label{ineq22}
B\ge(\log H_m)\exp\{4(m+1)(7+3\log{(m+1)})\},
\end{equation}
and
\begin{equation}\label{ineq23}
7+3\log(m+1)\ge\log k,
\end{equation}
then
\begin{equation*}
\abs{\Lambda}\ge\exp\{-c_2(m)k^{m+2}h_1\cdots h_m\log(\frac{2mB}{h_m})\},
\end{equation*}
where $c_2(m)$ is the function defined in Theorem \ref{thm2}.
\end{prop}

\section{Proof of Theorem \ref{thm2}}
In this section, we shall prove Theorem \ref{thm2}.  The proof is a standard
application of lower bounds for linear forms in logarithms of some special form.

Let $\OK=\BQ(\sqrt{D})$ and $\OO$ be the ring of integers.
We use the overline symbol to express the conjugate in $\OK$.
We denote by $\epsilon$ the fundamental unit in $\OK$ if $D>0$.

By the assumption that if $\gp$ is a prime ideal dividing both $[(X+Y\sqrt{D})/2]$
and $[(X-Y\sqrt{D})/2]$, then $\gp$ divides $[\sqrt{D}]$,
there exist some ideal factorizations $(m_i)=\gm_i \overline{\gm_i}$
for $i=1, \cdots, s$ and $A=\ga\overline{\ga}$ such that $[(X-Y\sqrt{D})/2]=\ga\gm_1^{e_1}\cdots\gm_s^{e_s}$.

Write $A^\prime=Am_1^{v_1}\cdots m_s^{v_s}$ and $e_i=hu_i+v_i$ with $0\leq v_i<h$.
Then we obtain  $[(X+Y\sqrt{D})/2]=(\alpha^\prime)(\mu_1)^{u_1}\cdots(\mu_s)^{u_s}$,
where $(\mu_i)=\gm_i^h$ for $1\leq i\leq s$ and $(\alpha^\prime)=\ga\gm_1^{v_1}\cdots\gm_s^{v_s}$,
each of which is necessarily principal.

{\bf Case 1.} $D>0$.  Denote by $\epsilon$ the fundamental unit of $\OK$ satisfying $\epsilon>1$.
By definition, $R=\log\epsilon$.  Then we can choose $\alpha^\prime$ and $\mu_i$ such that
$0<\overline{\alpha^\prime}\leq\alpha^\prime<\abs{\epsilon} A^{\prime 1/2}$
and $0<\overline{\mu_i}\leq\mu_i<\abs{\epsilon}m_i^{h/2}$.

Now there exists an integer $u_0$ such that
\begin{equation}
(X+Y\sqrt{D})/2=\alpha^\prime\epsilon^{u_0}\mu_1^{u_1}\cdots\mu_s^{u_s}
\end{equation}
and
\begin{equation}
(X-Y\sqrt{D})/2=\abs{\overline{\alpha^\prime}\epsilon^{-u_0}\overline{\mu_1}^{u_1}\cdots\overline{\mu_s}^{u_s}}.
\end{equation}

We assume that $u_0<0$ and put $b_0=-u_0$.  Clearly we have
\begin{equation}
Y=(\alpha^\prime\epsilon^{-b_0}\mu_1^{u_1}\cdots\mu_s^{u_s}-\abs{\overline{\alpha^\prime}\epsilon^{b_0}\overline{\mu_1}^{u_1}\cdots\overline{\mu_s}^{u_s}})/\sqrt{D}.
\end{equation}

Since $X$ and $Y$ are positive, we have $\abs{X-Y\sqrt{D}}<\abs{X+Y\sqrt{D}}$ and therefore
\begin{equation}
\abs{\overline{\alpha^\prime}\epsilon^{b_0}\overline{\mu_1}^{u_1}\cdots\overline{\mu_s}^{u_s}}<\abs{\alpha^\prime\epsilon^{-b_0}\mu_1^{u_1}\cdots\mu_s^{u_s}}.
\end{equation}
This clearly yields that
\begin{equation}\label{eq2a}
\abs{\alpha^\prime\epsilon^{-b_0}\mu_1^{u_1}\cdots\mu_s^{u_s}}>(Am_1^{e_1}m_2^{e_2}\cdots m_s^{e_s})^{1/2}.
\end{equation}
Thus we obtain
\begin{equation}\label{eq2b}
\abs{\Lambda}\leq\abs{Y/(\alpha^\prime\epsilon^{-b_0}\mu_1^{u_1}\cdots\mu_s^{u_s})}\sqrt{D}\leq c_0(Am_1^{e_1}m_2^{e_2}\cdots m_s^{e_s})^{-c_1}\sqrt{D},
\end{equation}
where
\begin{equation}
\Lambda=(\pm\frac{\overline{\alpha^\prime}}{\alpha^\prime})\epsilon^{2b_0}(\frac{\overline{\mu_1}}{\mu_1})^{u_1}\cdots(\frac{\overline{\mu_s}}{\mu_s})^{u_s}-1.
\end{equation}

By (\ref{eq2a}) and the choices of $\alpha'$ and $\mu_i$, we have
$h(\frac{\overline{\alpha^\prime}}{\alpha^\prime})\leq \frac{1}{2}\log A^\prime +R$, $h(\frac{\overline{\mu_i}}{\mu_i})\leq \frac{1}{2}h\log m_i +R$
and
\begin{equation}
\epsilon^{b_0}<\abs{\alpha^\prime\mu_1^{u_1}\cdots\mu_s^{u_s}}(Am_1^{e_1}m_2^{e_2}\cdots m_s^{e_s})^{-1/2}<\epsilon^{1+u_1+\cdots+u_s}.
\end{equation}
Hence $b_0\leq (u_1+\cdots+u_s)\leq su_j$, where $j$ is an index such that $u_j=\max u_i$.
Since $D\ge 5$, we have $m_i\ge 11$ and therefore $\log m_i\ge 2$.
From these estimates, we can apply Proposition \ref{prop1} with
\begin{equation*}
m=s+2,
\end{equation*}
\begin{equation*}
B=2su_j,
\end{equation*}
\begin{equation*}
h_i=\frac{1}{2}\log m_i +R,\mbox{ for }i=1,\cdots,s,
\end{equation*}
\begin{equation*}
h_{s+1}=\frac{1}{2}\langle R\rangle, 
\end{equation*}
\begin{equation*}
h_{s+2}=\frac{1}{2}\langle\log A^\prime +R\rangle,
\end{equation*}
to obtain
\begin{equation}\label{eq2c}
\abs{\Lambda}\ge\exp \{-c_1(s+2)2^{s+4}h_1\cdots h_{s+2}\log(\frac{2(s+2)B}{h_{s+2}})\}.
\end{equation}

Comparing (\ref{eq2b}) and (\ref{eq2c}), we obtain
\begin{equation}
(\frac{hc_1}{2s}\log m_j)B+\log(c_0^{-1}\abs{D}^{-\frac{1}{2}}A^{c_1})\leq c_1(s+2)2^{s+4}h_1\cdots h_{s+2}\log(\frac{2(s+2)B}{h_{s+2}}).
\end{equation}
Hence we have either
\begin{equation}
B\leq \frac{2s}{hc_1}\frac{\log(c_0\abs{D}^{\frac{1}{2}}A^{-c_1})}{\log m_j}
\end{equation}
or
\begin{equation}
(\frac{hc_1}{2s}\log m_j)B\leq c_2(s+2)2^{s+5}h_1\cdots h_{s+2}\log(\frac{2(s+2)B}{h_{s+2}}).
\end{equation}
In the former case, the inequality (\ref{ineq13a}) clearly holds.  Hence we limit to the latter case.
Multiplying both sides by $4s(s+2)(c_1hh_{s+2}\log m_j)^{-1}$, we have
\begin{equation}
\frac{2(s+2)B}{h_{s+2}}\leq \frac{s(s+2)c_2(s+2)2^{s+7}}{hc_1}h_1\cdots h_{s+1}(\log m_j)^{-1}\log(\frac{2(s+2)B}{h_{s+2}}).
\end{equation}
The assumption (\ref{ineq22}) gives
\begin{equation}
\frac{2(s+2)B}{h_{s+2}}\leq 2c_3\log c_3,
\end{equation}
where
\begin{equation}
c_3=\frac{s(s+2)c_2(s+2)2^{s+7}}{hc_1}h_1\cdots h_{s+1}(\log m_j)^{-1}.
\end{equation}
Hence we obtain
\begin{equation}
2su_j=B\leq c_4\langle\log A^\prime+2R\rangle\log\{(s+2)c_4\},
\end{equation}
where
\begin{equation}
c_4=\frac{2^5 s c_2(s+2)}{hc_1}\langle R\rangle(1+\frac{2R}{\log m_j})\prod_{i\ne j}(\log m_i+2R).
\end{equation}
Noting that $\max e_i\leq hu_i+h-1$, this yields the inequality (\ref{ineq13a}).

The case $u_0\geq 0$ remains to consider.  Let
\begin{equation}
\abs{\Lambda^\prime}=(\pm\frac{\alpha^\prime}{\overline{\alpha^\prime}})\epsilon^{2u_0}(\frac{\mu_1}{\overline{\mu_1}})^{u_1}\cdots(\frac{\mu_s}{\overline{\mu_s}})^{u_s}-1.
\end{equation}
Then we easily see that by the hypothesis
\begin{equation}\label{eq221}
\abs{\Lambda^\prime}\leq \frac{2c_0(Am_1^{e_1}m_2^{e_2}\cdots m_s^{e_s})^{-c_1}\sqrt{D}}{1-c_0(Am_1^{e_1}m_2^{e_2}\cdots m_s^{e_s})^{-c_1}\sqrt{D}}
\end{equation}

We easily see that $\mu_i\geq \overline{\mu_i}$ and $\mu_i/\overline{\mu_i}\geq (t_i+\sqrt{D})/(t_i-\sqrt{D})$ with $t_i=2m_i^{h/2}$.  Hence we obtain
\begin{equation}
\abs{\Lambda}\geq (\frac{t_s+\sqrt{D}}{t_s-\sqrt{D}})^E-1,
\end{equation}
where $E=\max u_i$ in $1\leq i\leq s$.  We have
\begin{equation}\label{eq222}
\abs{\Lambda}\geq \frac{2E\sqrt{D}}{t_s}.
\end{equation}

From (\ref{eq221}) and (\ref{eq222}) we have
\begin{equation}
c_0(Am_1^{e_1}m_2^{e_2}\cdots m_s^{e_s})^{-c_1}\sqrt{D}\geq \min\{1/2, \frac{E\sqrt{D}}{2t_s}\}.
\end{equation}
This implies either
\begin{equation}\label{eq223}
Am_1^{e_1}m_2^{e_2}\cdots m_s^{e_s}\leq (2c_0\sqrt{D})^{1/c_1}
\end{equation}
or
\begin{equation}\label{eq224}
AEm_1^{e_1}m_2^{e_2}\cdots m_s^{e_s}\leq (2c_0t_s)^{1/c_1}.
\end{equation}

If (\ref{eq223}) holds, then the estimate $m_i\geq 11\geq e^2$ gives the inequality (\ref{ineq13a}).  If (\ref{eq224}) holds, then we have
\begin{equation}
E\leq\frac{h\log(2c_0m_s^{1/2})}{\log m_1},
\end{equation}
which immidiately yields the inequality (\ref{ineq13a}).  This completes the proof in the case $D>0$.

{\bf Case 2.} $D<0$. In this case, we have $\abs{\mu_i}=\abs{\overline{\mu_i}}=m_i^{h/2}$ and
$\abs{\alpha^\prime}=\abs{\overline{\alpha^\prime}}=A^\prime$.

Now there exists an integer $u_0$ such that
\begin{equation}
(X+Y\sqrt{D})/2=\alpha^\prime\omega^{u_0}\mu_1^{u_1}\cdots\mu_s^{u_s}
\end{equation}
and
\begin{equation}
(X-Y\sqrt{D})/2=\overline{\alpha^\prime}\omega^{-u_0}\overline{\mu_1}^{u_1}\cdots\overline{\mu_s}^{u_s},
\end{equation}
where $\omega$ is a primitive sixth root of unity when $q=3$ and is $-1$ otherwise.  Note that $\omega\overline{\omega}=1$.

\begin{equation}
Y=(\alpha^\prime\omega^{u_0}\mu_1^{u_1}\cdots\mu_s^{u_s}-\overline{\alpha^\prime}\omega^{-u_0}\overline{\mu_1}^{u_1}\cdots\overline{\mu_s}^{u_s})/\sqrt{D}.
\end{equation}

It is clear that
\begin{equation}
\abs{\overline{\alpha^\prime}\omega^{-u_0}\overline{\mu_1}^{u_1}\cdots\overline{\mu_s}^{u_s}}=\abs{\alpha^\prime\omega^{u_0}\mu_1^{u_1}\cdots\mu_s^{u_s}}=(Am_1^{e_1}m_2^{e_2}\cdots m_s^{e_s})^{1/2}.
\end{equation}
Hence we obtain
\begin{equation}\label{eq2d}
\abs{\Lambda}\leq\abs{Y/(\alpha^\prime\omega^{u_0}\mu_1^{u_1}\cdots\mu_s^{u_s})}\sqrt{\abs{D}}\leq c_0(Am_1^{e_1}m_2^{e_2}\cdots m_s^{e_s})^{-c_1}\sqrt{\abs{D}},
\end{equation}
where
\begin{equation}
\Lambda=\omega^{-2u_0}(\frac{\overline{\alpha^\prime}}{\alpha^\prime})(\frac{\overline{\mu_1}}{\mu_1})^{u_1}\cdots(\frac{\overline{\mu_s}}{\mu_s})^{u_s}-1.
\end{equation}

By (\ref{eq2a}) and the choices of $\alpha'$ and $\mu_i$, we have
$h(\frac{\overline{\alpha^\prime}}{\alpha^\prime})=\frac{1}{2}\log A^\prime$, $h(\frac{\overline{\mu_i}}{\mu_i})\leq \frac{1}{2}h\log m_i$
and
From these estimates, we can apply Proposition \ref{prop1} with
\begin{equation*}
m=s+1,
\end{equation*}
\begin{equation*}
B=u_j,
\end{equation*}
\begin{equation*}
h_i=\frac{1}{2}\langle\log m_i\rangle,\mbox{ for }i=1,\cdots,s,
\end{equation*}
\begin{equation*}
h_{s+1}=\frac{1}{2}\langle\log A^\prime\rangle,
\end{equation*}
to obtain
\begin{equation}\label{eq2e}
\abs{\Lambda}\ge\exp \{-c_2(s+1)2^{s+3}h_1\cdots h_{s+1}\log(\frac{2(s+1)B}{h_{s+1}})\}.
\end{equation}
Comparing (\ref{eq2d}) and (\ref{eq2e}), we obtain
\begin{equation}
(hc_1\log m_j)B+\log(c_0^{-1}\abs{D}^{-\frac{1}{2}}A^{c_1})\leq c_2(s+1)2^{s+3}h_1\cdots h_{s+1}\log(\frac{2(s+1)B}{h_{s+1}}).
\end{equation}
Thus we have either
\begin{equation}
B\leq \frac{\log(c_0\abs{D}^{\frac{1}{2}}A^{-c_1})}{hc_1\log m_j}
\end{equation}
or
\begin{equation}
(hc_1\log m_j)B\leq c_2(s+1)2^{s+4}h_1\cdots h_{s+1}\log(\frac{2(s+1)B}{h_{s+1}}).
\end{equation}
In the former case, the inequality (\ref{ineq13b}) clearly holds.  Hence we limit to the latter case.
Multiplying both sides by $4s(s+1)(c_1hh_{s+1}\log m_j)^{-1}$, we have
\begin{equation}
\frac{2(s+1)B}{h_{s+1}}\leq \frac{(s+1)c_2(s+1)2^{s+5}}{hc_1}h_1\cdots h_s(\log m_j)^{-1}\log(\frac{2(s+1)B}{h_{s+1}}).
\end{equation}
The assumption (\ref{ineq22}) gives
\begin{equation}
\frac{2(s+1)B}{h_{s+1}}\leq 2c_5\log c_5,
\end{equation}
where
\begin{equation}
c_5=\frac{(s+1)c_2(s+1)2^{s+5}}{hc_1}h_1\cdots h_s(\log m_j)^{-1}.
\end{equation}
Since $h_1/\log m_j\leq h_1/\log m_1\leq 2/\log 7$ and $\log m_i\ge \log 13>2$ for $i>1$, we obtain
\begin{equation}
B\leq c_6(\log A^\prime)\log\{(s+1)c_6\},
\end{equation}
where
\begin{equation}
c_6=\frac{2^6 c_2(s+1)}{hc_1\log 7}\prod_{i>1}(\log m_i).
\end{equation}
Noting that $\max e_i\leq hB+h-1$, this yields the inequality (\ref{ineq13b}).  This completes the proof.

\section{Proof of the Theorem \ref{thm1}}
We put $D=(-1)^{(q-1)/2}q$.  It is clear that $D\equiv 1\pmod{4}$.
As in the previous section, we denote by $\OK$ and $\OO$, respectively, $\BQ(\sqrt{D})$ and its ring of integers.
We use the overline symbol to express the conjugate in $\OK$.
Assume $x\ge q^{3/2}$.
Let
\begin{equation}
P^+(x)=\prod_{(\frac{m}{q})=1} (x-\zeta^m)=\sum_{i=0}^{\frac{q-1}{2}}a_i x^{\frac{q-1}{2}-i}
\end{equation}
and
\begin{equation}
P^-(x)=\prod_{(\frac{m}{q})=-1} (x-\zeta^m),
\end{equation}
where $\zeta$ is a primitive $q$-th root of unity.
Then it is well-known that $P^+(x)=\overline{P^-}(x)$ has its coefficients in $\OO$.
We have $\abs{a_i}\leq\left(\begin{array}{cc}(q-1)/2\\i\end{array}\right)\leq q^i$.

We put
\begin{equation*}
X=f(x)=(P^+(x)+P^-(x))
\end{equation*}
and
\begin{equation*}
Y=g(x)=(P^+(x)-P^-(x))/(\sqrt{D}).
\end{equation*}
Then it is clear that the coefficients of $f(x)$ and $g(x)$ belong in $\BZ$ and $g(x)$ has degree $(q-3)/2$.
So we can write
\begin{equation}
g(x)=\sum_{i=0}^{\frac{q-3}{2}}b_i x^{\frac{q-3}{2}-i}.
\end{equation}
By a well-known result on Gaussian sums, we observe $b_0=\pm 1$.
Moreover, we have $\abs{b_i}=\abs{a_{i+1}}/(2\sqrt{\abs{D}})\leq q^{i+1/2}/2$.
Hence we have
\begin{equation}
\abs{g(x)}\leq x^{(q-3)/2}+q^{3/2}x^{(q-5)/2}\leq 2x^{(q-3)/2}<2(\frac{x^q-1}{x-1})^{\frac{q-3}{2(q-1)}}.
\end{equation}
Moreover, we observe that if a prime ideal $\gp$ divides $[P^+(x), P^-(x)]$, then
$\gp$ divides $[\sqrt{D}]$, for $\gp$ divides $[x-\zeta^i, x-\zeta^j]$
for some $i$, $j$ with $i\ne j$ and therefore $\gp$ must divide $[\zeta^i-\zeta^j]$.
Based on these facts and the identity
\begin{equation}
4\frac{x^q-1}{x-1}=X^2-DY^2,
\end{equation}
we can apply Theorem \ref{thm2} with $(c_0, c_1)=(2, 1/q)$ and we immediately
obtain the inequalities (\ref{ineq11a}) and (\ref{ineq11b}).

\section{A simple estimate for $U$}

In this section, we give a simple estimate for $U$ since the definition of $U$ in Theorem \ref{thm1} is too complicated for application.  We begin by noting that $h, R\leq q^{1/2}\log 4q$(\cite{Fai}).

{\bf Case 1.} $q\equiv 1\pmod{4}$.

We have
\begin{equation}
U<\max\{\frac{h^2q\log(4m_s^{1/2})}{\log m_1}, \frac{f\log 2q^{1/2}}{\log m_1}, C_1\langle\log A^\prime\rangle(\log\frac{(s+1)C_1}{h})\}+h-1.
\end{equation}

Firstly, if
\begin{equation}
U<\frac{h^2q\log(4m_s^{1/2})}{\log m_1}+h-1,
\end{equation}
then we have $U<q^5\log m_s$ if $s>1$ and $U<q^5$ if $s=1$.  Secondly, if
\begin{equation}
U<\frac{q\log 2q^{1/2}}{\log m_1}+h-1,
\end{equation}
then we have $U<q^5$.  Finally we consider the case
\begin{equation}
U<C_1\langle\log A^\prime\rangle(\log\frac{(s+1)C_1}{h})+h-1.
\end{equation}
We have
\begin{equation}
\langle\log A^\prime\rangle\leq \max\{2, (h-1)\log(m_1\cdots m_s)\}\leq q^{1/2}\log(4q)\log(m_1\cdots m_s),
\end{equation}
\begin{equation}
\log[(s+1)C_1]\leq 2\log C_1\leq C_1.
\end{equation}
Thus
\begin{eqnarray}
\log U&\leq&\log[C_1\langle\log A^\prime\rangle(\log\frac{(s+1)C_1}{h})+h-1]\\
&\leq&\log[C_1^2 q^{1/2}\log(4q)\log(m_1\cdots m_s)+q^{1/2}\log(4q)]\\
&\leq&\log[(\left(\frac{2^6qc_2(s+1)}{\log 7}\prod_{i=2}^{s}(\log m_i)\right)^2\log(m_1\cdots m_s)+1) q^{1/2}\log(4q)]\\
&\leq& 2\log[2^6qc_2(s+1)\prod_{i=2}^{s}(\log m_i) q^{1/2}\log(4q)\log(m_1\cdots m_s)]\\
&\leq& 12\log2+2\log 1500+2(s+2)\log 38+(6s+4)\log(s+2)\\
&+&2\sum_{i=2}^{s}\log\log m_i+\log\log(m_1\cdots m_s)+3\log q+2\log\log(4q)\\
&\leq&\log c_7+2s\log 38+(6s+4)\log(s+2)+3\sum_{i=2}^{s}\log\log m_i+5\log q\\
&\leq&\log c_7+(13s+4)\log(s+2)+3\sum_{i=2}^{s}\log\log m_i+5\log q.
\end{eqnarray}

{\bf Case 2.} $q\equiv 3\pmod{4}$.
We have
\begin{equation}
U<\max\{\frac{q\log 2q^{1/2}}{\log m_1}, C_0\langle\log A^\prime+2R\rangle(\log\frac{2s(s+2)C_0}{h})\}+h-1,
\end{equation}
If
\begin{equation}
U<\frac{q\log 2q^{1/2}}{\log m_1}+h-1,
\end{equation}
then it follows that $U<q^3$.  We next consider the case
\begin{equation}
U<C_0\langle\log A^\prime+2R\rangle(\log\frac{2s(s+2)C_0}{h})+h-1.
\end{equation}
We have
\begin{equation}
\langle\log A^\prime\rangle+2R\leq q^{1/2}\log(4q)[2+\log(m_1\cdots m_s)],
\end{equation}
\begin{equation}
\log[2s(s+2)C_0]\leq 2\log C_0\leq C_0.
\end{equation}
Thus
\begin{eqnarray}
\log U&\leq&\log[C_0\langle\log A^\prime+2R\rangle(\log\frac{s(s+2)C_0}{h})+h-1]\\
&\leq&\log[C_0^2 q^{1/2}\log(4q)(2+\log[m_1\cdots m_s])+q^{1/2}\log(4q)]\\
&\leq&\log[(C_0^2\log[m_1\cdots m_s]+3)q^{1/2}\log(4q)]\\
&\leq& 2\log[2^5qc_2(s+2)\prod_{i=2}^{s}(\log m_i) q^{1/2}\log(4q)\log(m_1\cdots m_s)]\\
&\leq& 12\log2+2\log 1500+2(s+2)\log 38+(6s+4)\log(s+2)\\
&+&2\sum_{i=2}^{s}\log\log m_i+\log\log(m_1\cdots m_s)+3\log q+2\log\log(4q)\\
&\leq&\log c_7+2s\log 38+(6s+4)\log(s+2)+3\sum_{i=2}^{s}\log\log m_i+5\log q\\
\end{eqnarray}

In both cases, we obtain
\begin{equation}\label{eq51}
\log U\leq\log c_7+2s\log 38+(6s+4)\log(s+2)+3\sum_{i=2}^{s}\log\log m_i+5\log q.
\end{equation}

\section{Proof of Theorem \ref{thm3}}

Let $(x^q-1)/(x-1)=p_1^{e_1}p_2^{e_2}\cdots p_s^{e_s}$, where we denote the prime divisors
of $(x^q-1)/(x-1)$ by $p_1<p_2<\cdots<p_s$.  By definition, we have $P[(x^q-1)/(x-1)]=p_s$.
Write for brevity, $P=p_s$.

Since $p_i=q$ or $p_i\equiv 1\pmod{q}$, we have
\begin{equation}
s\le\frac{P}{(q-1)\log P}+c\frac{P}{(\log P)^2}
\end{equation}
for some effectively computable constant $c$ depending on $q$(\cite[Theorem 9.6]{Kar}).
Moreover, we have an trivial estimate $(x^q-1)/(x-1)\le P^{sU}$ and therefore
$x\le P^{sU/(q-1)}$.  Now the inequality (\ref{eq51}) immidiately gives
\begin{equation}\label{eq61}
P>(\frac{q-1}{6}-\epsilon)\log \log x
\end{equation}
for $P>P_0$, where $P_0$ denotes an effectively computable constant depending only on $q$ and
$\epsilon$.  We see that this fact implies that (\ref{eq61}) also holds for $x>x_0$,
since $P$ tends to infinity together with $x$.  This proves the theorem.

\section{Proof of Theorem \ref{thm4}}

In order to prove Theorem \ref{thm4}, we prove a combinatorial lemma concerning
the distribution of the solutions of (\ref{eq14}).

\begin{lem}\label{lm41}
Let $p_0$, $p_1$, $p_2$ be distinct primes, $e$ and $q$ be a positive integers.
Put $H_i=e\log p_0/\log p_i$ for $i=1, 2$.
If the equation
\begin{equation}\label{eq71}
p_i^q\equiv 1\pmod{p_0^e}
\end{equation}
holds for $i=1, 2$, then
\begin{equation}
\frac{3}{4}H_1H_2\leq (q, p_0-1).
\end{equation}
\end{lem}
\begin{proof}
Consider the congruence
\begin{equation}\label{eq73}
p_1^{a_1}p_2^{a_2}\equiv 1\pmod{p_0^e}
\end{equation}
with $0\leq a_i\leq H_i$.
We assume that (\ref{eq73}) has no solution.  Then $p_1^{a_1}p_2^{a_2}$
takes distinct values $\pmod{q^e}$ for each $a_1$ and $a_2$.
But (\ref{eq71}) implies that this takes at most $(q, p_0-1)$ distinct values.
Hence we obtain $(\floor{H_1}+1)(\floor{H_2}+1)\leq (q, p_0-1)$.

The remaining case is when (\ref{eq73}) has a solution $(a_1, a_2)$.
Consider the congruence
\begin{equation}\label{eq74}
p_1^{b_1}p_2^{b_2}\equiv p_1^{c_1}p_2^{c_2}\pmod{p_0^e}
\end{equation}
with $0\leq b_i, c_i\leq H_i$ and $(b_1, b_2)\ne(c_1, c_2)$.
If (\ref{eq74}) has a solution, then we have $p_1^{c_1-b_1}\equiv p_2^{b_2-c_2}\pmod{p_0^e}$.
Hence
\begin{equation}\label{eq75}
c_i\leq b_i\mbox{ for both }i=1, 2\mbox{ or }c_i\ge b_i\mbox{ for both }i=1, 2.
\end{equation}
In both case we have
\begin{equation}
p_1^{\abs{c_1-b_1}}p_2^{\abs{c_2-b_2}}\equiv p_1^{a_1}p_2^{a_2}\equiv 1\pmod{p_0^e}.
\end{equation}
We shall show that $\abs{c_i-b_i}\ge a_i$ for both $i=1, 2$.  Otherwise we have $a_1>\abs{c_1-b_1}>0$
and $a_2>\abs{c_2-b_2}>0$ by virtue of (\ref{eq75}).  Hence $a_1\log p_1 +a_2\log p_2>2p_0^e$,
which is imcompatible with the ranges of $a_1$ and $a_2$.
Thus we conclude that $p_1^{b_1}p_2^{b_2}$ takes distinct values for each $b_1, b_2$
satisfying $0\leq b_i\leq H_i$ for $i=1, 2$ and $b_i<a_i$ for at least one $i$.
Hence we obtain $\frac{3}{4}H_1H_2\leq (q, p_0-1)$.  This completes the proof.
\end{proof}
\begin{lem}\label{lm42}
Let $p_0$, $p_1$, $p_2$ be distinct primes with $p_2>p_1$ and $q$ be a positive integer.
If there are integers $e_i$ such that $p_0^{e_i}\mid\sigma(p_i^{q-1})$ and $p^{e_i}\ge\sigma(p_i^{f-1})^{1/s}$ for $i=1, 2$, then
\begin{equation}\label{eq76}
\log p_2>\frac{3(q-1)^2}{4qs^2}\log p_1.
\end{equation}
\end{lem}
\begin{proof}
Let $e=\min\{e_1, e_2\}$ and $H_i=e\log p_0/\log p_i$.  Then it is clear that $H_i\ge(q-1)/s$.
By Lemma \ref{lm41}, we obtain
\begin{equation}
\log p_2\ge\frac{3}{4q}H_1^2\log p_1\ge \frac{3(q-1)^2}{4qs^2}\log p_1.
\end{equation}
This proves the lemma.
\end{proof}
If $q$ is sufficiently large compared to $s$, then we obtain a more simple inequality.
\begin{cor}\label{cr41}
If the conditions in Lemma \ref{lm42} hold and $q>\frac{16}{9}es^4$, then
\begin{equation}\label{eq77}
\log p_2>q^{1/2}\log p_1.
\end{equation}
\end{cor}

If the equation (\ref{eq14}) holds, then there exists an index $i$ such that $m_i^{e_i}\ge\sigma(r^{q-1})^{1/s}$.
So we can divide the solutions of the equation (\ref{eq14}) into $s$ sets $R_1, \cdots, R_s$ so that
if $(r, e_1, \cdots, e_s)$ is the solution of (\ref{eq14}) with $r\in R_i$, then $m_i^{e_i}\ge\sigma(r^{q-1})^{1/s}$.
By Corollary \ref{cr41}, if $r_1<r_2<\cdots$ are the elements of $R_i$, then $\log r_{j+1}>f^{1/2}\log r_j$.

Now we have $\sigma(r_1^{q-1})\ge m_i$ and $\sigma(r_{\# R_i}^{q-1})\leq m_i^{sU}$.
It immidiately follows that $r_{\# R_i}<r_1^{sqU/(q-1)}$.  These facts give
\begin{equation}
\#R_i\leq\frac{2\log sqU/(q-1)}{\log q}+1
\end{equation}
and therefore the number of solutions of the equation (\ref{eq14}) is at most
\begin{equation}\label{eq78}
s(\frac{2\log sqU/(q-1)}{\log q}+1).
\end{equation}

Applying the upper bound for $U$ given in (\ref{eq51}), we obtain

\begin{eqnarray}
&s&(\frac{2\log s+\log q/(q-1)+\log U}{\log q}+1)\\
\leq&s&(\frac{\log c_7+(13s+6)\log(s+2)+3\sum_{i=2}^{s}\log\log m_i
+6\log q}{\log q}+1)\\
\leq&s&(\frac{\log c_7+19s\log(s+2)+3\sum_{i=2}^{s}\log\log m_i}{\log q}+7).
\end{eqnarray}

This completes the proof.

{}
\vskip 12pt

{\small Tomohiro Yamada}\\
{\small Department of Mathematics,}\\
{\small Graduate School of Science,}\\
{\small Kyoto University, Kyoto, 606-8502, Japan}\\
{\small e-mail: \protect\normalfont\ttfamily{tyamada@math.kyoto-u.ac.jp}}

\begin{thebibliography}{}
\bibitem{BG1}
Y. Bugeaud and K. Gy\H{o}ry.
Bounds for the solutions of unit equations.
\emph{Acta Arith.}\textbf{74} (1996), 67--80.
\bibitem{BG2}
Y. Bugeaud and K. Gy\H{o}ry.
Bounds for the solutions of Thue-Mahler equations and norm form equations.
\emph{Acta Arith.}
\textbf{74} (1996), 273--292.
\bibitem{BHM}
Y. Bugeaud, G. Hanrot and M. Mignotte.
Sur l'{\'e}quation diophantienne $\frac{x^n-1}{x-1}=y^m$, III.
\emph{Proc London Math. Soc.}
\textbf{84} (2002), 59--78.
\bibitem{Fai}
A. Faisant.
\emph{L'{\'e}quation diophantienne du second degr{\'e}}
(Herrmann, Paris, 1991).
\bibitem{Gyo}
K. Gy\H{o}ry.
On some arithmetical properties of Lucas and Lehmer numbers.
\emph{Acta Arith.}
\textbf{40} (1982), 369--373.
\bibitem{Kar}
A. A. Karatsuba.
\emph{Basic Analytic Number Theory}
(Springer-Verlag, 1993).
\bibitem{Kot}
S. V. Kotov.
The Thue-Mahler equation in relative fields(Russian),
\emph{Acta Arith.}
\textbf{27}(1975), 293--315.
\bibitem{Wal}
M. Waldschmidt.
Minorations de combinaisons lin{\'e}aires de logarithmes de nombres alg{\'e}briques.
\emph{Canadian J. Math.}
\textbf{45} (1993), 176--224.
\end{thebibliography}
\end{document}